\documentclass[11pt,a4paper]{article}
\usepackage[latin1]{inputenc} 
\usepackage[english]{babel}
\usepackage[T1]{fontenc} 
\usepackage{amsfonts,amssymb,latexsym,epsfig,amsmath}
\usepackage{graphicx}
\usepackage{amsthm}
\usepackage{dsfont}
\usepackage{float}
\usepackage[font=small,labelfont=bf]{caption}
\usepackage[font=small,labelfont=]{subfig}
\usepackage{color}
\usepackage{hyperref}

\newtheorem{theorem}{Theorem}[section]
\newtheorem{proposition}[theorem]{Proposition}
\newtheorem{corollary}[theorem]{Corollary}
\newtheorem{lemma}[theorem]{Lemma} 

\newtheorem{remark}[theorem]{Remark}
\numberwithin{equation}{section}

\newcommand{\Beq}{\begin{equation}}
\newcommand{\Eeq}{\end{equation}}
\newcommand{\BS}{\begin{subequations}}
\newcommand{\ES}{\end{subequations}}
\newcommand{\Beqn}{\begin{equation*}}
\newcommand{\Eeqn}{\end{equation*}}
\newcommand{\Beqa}{\begin{eqnarray}}
\newcommand{\Eeqa}{\end{eqnarray}}
\newcommand{\Beqan}{\begin{eqnarray*}}
\newcommand{\Eeqan}{\end{eqnarray*}}
\newcommand{\N}{\hbox{$ I\kern -0.23em N$}}

\newcommand{\R}{{\mathbb R}}

\newcommand{\eps}{\varepsilon}

\newcommand{\md}{\mathfrak d}

\begin{document}
\title
{Steady periodic water waves bifurcating for fixed-depth rotational flows with discontinuous vorticity}

\author{David Henry and Silvia Sastre-Gomez \\[.5em]
\small School of Mathematical Sciences, \\ \small University College Cork, \\ \small Cork, Ireland.}

\date{}

\maketitle

\begin{abstract}
\noindent In this article we apply local bifurcation theory to prove the existence of small-amplitude steady periodic water waves, which propagate over a flat bed with a specified fixed mean-depth, and where the underlying flow has a discontinuous vorticity distribution. 
\end{abstract}
%
\section{Introduction}
In this paper we prove the existence of small-amplitude steady periodic water waves which propagate over a flat bed with a specified fixed mean-depth, and where the underlying flow has a discontinuous vorticity distribution. This existence result is achieved by way of an implementation of the Crandall-Rabinowitz local bifurcation theory applied to a generalised version of a modified-height function formulation of the governing equations. This novel reformulation, introduced in \cite{Hen13,DH1}, enables us to maintain the mean fluid depth as a fixed parameter in the bifurcation analysis, as opposed to the (by now) standard approach initiated in  \cite{CS2} whereby, following the Dubreil-Jacotin \cite{DubJac} transformation, the mean depth does not feature in the bifurcation setting. Instead, in \cite{CS2,CS11} the mass-flux is maintained as a fixed-parameter in the standard height function formulation. As in many complex bifurcation problems the role of the involved parameters is by no means canonical, and this is particularly true for the water wave problem. Fixing the mean depth of the wave is heuristically and physically quite natural, and furthermore since numerical studies \cite{KoS1} suggest that for fixed mass-flux the mean depth is not a constant of motion, the present contribution represents in a sense a completion of the analysis of the mathematical problem.

The seminal work \cite{CS2} was the first mathematically rigorous work which established the existence of water waves with general vorticity distributions, prior to this much of the rigorous mathematical analysis of water waves focussed on irrotational flows--- cf. \cite{Tol} for a nice survey of the rigorous mathematical analysis of irrotational waves, while surveys of recent progress in the mathematical analysis of rotational flows can be found in \cite{Cons_book,ConJPA}. 
Following \cite{CS2}, mathematically rigorous existence theory for rotational flows underlying steady periodic water waves has thrived, with a proliferation of results modelling scenarios with increasing levels of physical complexity, such as incorporating fluid stratification, the capillary effects of surface tension, the presence of weak stagnation and critical points, and  even accommodating the possibility of overhanging wave-profiles, cf. \cite{CS11,ConVar,EMM1,HM14x,HM13,MarMat,M14x,MMexist,Wa6b,Wah3,Wa1,Wa3,Wa4}. The practical relevance of flows with vorticity has long been recognised, particularly as they  serve as models for wave-current interactions, among other intricate and physically important phenomena \cite{Cons_book,Jon,ThomKlop}. The extension of the existence result in \cite{DH1} to allow for discontinuous vorticity is both mathematically interesting and challenging. Physically, it is known that wave-current interaction is most dramatically observed where there is a rapid change in the current strength \cite{Jon}, a scenario which may be crudely modelled by  a discontinuous vorticity distribution. Examples of such processes include wind-generated waves, and near-bed  current layers which  account  for sediment transport. 

This paper is laid out as follows. In Section~\ref{SecPre} we recast the governing equations in terms of a generalised version of the modified height function formulation, in the process establishing the equivalence of this new formulation with the weak versions of the standard Euler equation and stream function formulations, and additionally we present strong {\em a priori} regularity results for the streamlines (including the unknown free-boundary) of the resulting flow. The main result of our paper, which establishes the existence of small-amplitude weak solutions, is  presented in Theorem ~\ref{MR}, and Section~\ref{SecLB} is then dedicated to proving this result. Finally, in Section~\ref{SecEx} we present explicit, physically-motivated examples of various flows for which our bifurcation result applies. 
 
\section{Preliminaries}\label{SecPre}
\subsection{Formulations of the governing equations}
We consider two-dimensional steady periodic travelling surface waves propagating over water of mean depth $d>0$  where the predominant external restoration force is gravity. The $(x,y)-$coordinate system is fixed by locating the  mean--water level at $y=0$, which implies that the mean of the waves' free surface over a wavelength must equal zero: if $\eta(x,t)$ denotes the unknown free surface at some given fixed-time $t$, then 
\Beq\label{eta}
\int \eta(x,t)dx=0,
\Eeq
where the integral is taken over a wavelength. The impermeable flat bed is given by $y=-d$. Aside from the physically--enforced restriction \eqref{eta}, the free surface $\eta$ is {\em a priori} undetermined and is therefore an unknown in the mathematical formulation of the water wave problem. 
We seek steady travelling waves propagating with a constant wave speed $c>0$ in the positive $x-$direction. Accordingly, all functions must have a $x,t$ relationship of the form  $x-ct$, and the change of coordinates $(x-ct,y)\mapsto (x,y)$ transforms to a steady reference frame moving alongside the wave where the flow is time independent. We denote the closure of the fluid domain by $\overline{D_{\eta}}=\{(x,y)\in \mathbb R^2: -d\leq y\leq \eta (x)\}$. The governing equations for the motion of the perfect (inviscid and incompressible) fluid take the form of the Euler equations together with boundary conditions given by 
\begin{subequations}\label{euler_n}
\begin{align}
\label{masscon} u_x+v_y&=0 & \mbox{ in } & D_{\eta}, \\
\label{E1} (u-c)u_x+vu_y&=-P_x  & \mbox{ in }  & D_{\eta}, \\  
\label{E2} (u-c)v_x+vv_y&=-P_y-g & \mbox{ in } &  D_{\eta}, \\
\label{ksc'} v&=(u-c)\eta_x & \mbox{ on } & y=\eta(x),   \\
\label{bnddyn1}  P&=P_{atm} &\mbox{ on } & y=\eta(x), \\
\label{kbc}v&=0  & \mbox{ on } & y=-d,   
\end{align}
where $(u,v)$ is the velocity field, $P(x,y)$ is the pressure distribution function, $P_{atm}$ is the constant atmospheric pressure and  $g$ is the gravitational constant of acceleration.
We now make the additional  assumption  that there are no weak-stagnation points, that is,  
\begin{equation}\label{umax} 
u<c
\end{equation}
\end{subequations}
throughout the fluid. This is a physically reasonable assumption for water waves, without underlying currents containing strong non-uniformities, and which are not near breaking \cite{J90,Light}. These flows do not contain any stagnation points, and the 
individual fluid particles move with a horizontal velocity which is less than the 
speed with which the surface wave propagates.  The non-stagnation condition is also  essential mathematically in our  reformulations of the water wave equations, as we see below. For two-dimensional motion the vorticity is given by
\begin{equation}
\label{vort}
\omega=u_y-v_x.
\end{equation} 
We work with periodic waves and note that  we may effectively choose  the period to be equal to $2\pi$ without loss of generality in subsequent considerations. This is due to the fact that if we are dealing with water waves of wavelength $L$ in the governing equations \eqref{masscon}--\eqref{kbc} then performing the scaling of variables 
\Beqn
(x,y,t,g,\omega,\eta,u,v,P,c) \mapsto (\kappa x,\kappa y, \kappa t,  \kappa ^{-1}g,  \kappa ^{-1}\omega, \kappa \eta, u,v,P,c), \Eeqn
where $\kappa =\frac{2\pi}{L}$ is the wavenumber, we end up with a $2\pi-$periodic system in the new variables identical to \eqref{euler_n} except $g,\omega$ are replaced by $\kappa ^{-1}g,\kappa ^{-1}\omega$. 
We can write the Euler equations \eqref{E1}-\eqref{E2} coupled with the mass conservation 
equation \eqref{masscon} and the boundary conditions \eqref{ksc'}-\eqref{kbc} 
in the (weak) divergence form as 
\begin{subequations}\label{euler_weak}
\begin{align}
	\label{weul1}
	-cu_x+\left(u^2\right)_x+\left(uv\right)_y&=-P_x   & \mbox{ in } & D_{\eta},\\ 
	\label{weul2}
	-cv_x+\left(uv\right)_x+\left(v^2\right)_y&=-P_y-g  & \mbox{ in } & D_{\eta},\\
	\label{weul3}
	u_x+v_y&=0 & \mbox{ in } & D_{\eta},\\
	\label{weul4}
	v&=0  & \mbox{ on } & y=-d,  \\
	\label{weul5} 
	v&=(u-c)\eta_x  & \mbox{ on } & y=\eta(x),  \\
	\label{weul6}
	P&=P_{atm}  & \mbox{ on } & y=\eta(x).
\end{align}
\end{subequations}   
With this weak formulation we may consider solutions with weaker 
regularity than those of \eqref{euler_n}. The equations \eqref{weul1}--\eqref
{weul3} will be understood in the sense of distributions, whereas the 
boundary conditions  \eqref{weul4}--\eqref{weul6} will be understood in the 
classical sense. The type of solutions of \eqref{euler_weak} we are 
interested in are  H\"older continuously 
differentiable functions $u,v,P\in W^{1,r}_{per}(D_{\eta})\subset C^{0,\alpha}_{per}(\overline{D_{\eta}})$, 
with $\eta\in C^{1,\alpha}_{per}(\R)$, for some H\"older exponent $\alpha\in (0,1)$ and 
$r=2/(1-\alpha)$, where $r$ is chosen to ensure the embedding $W^{1,r}
_{per}(D_{\eta})\subset C^{0,\alpha}_{per}(\overline{D_{\eta}})$. The \emph{per} subscript indicates that our solutions are $2\pi-$periodic, and furthermore we assume that solutions have a single crest located at $x=0$ and 
troughs located at $x=\pm \pi$, and the condition \eqref{umax} on $u$ and $c
$ will hold throughout the fluid domain. 

The governing equations \eqref{euler_n} have a second and, as we mention below, equivalent formulation  when recast in terms of the  stream-function. The stream function $\psi$ is defined up to a constant by the relations 
\begin{equation}\label{psi1}
	\psi_y=u-c,\quad \psi_x=-v.
\end{equation}
The constant may be fixed by setting $\psi=0$ on $y=\eta(x)$, and we observe by expressing
\[
	\psi(x,y)=-p_0+\int_{-d}^{y}(u(x,s)-c)ds,
\]
 that $\psi$ is also periodic with period $2\pi$. We note that one consequence of the definition \eqref{psi1} is that the level sets of the stream 
function $\psi(x,y)$ are streamlines of the fluid motion, hence the name. 
The  condition \eqref{umax} ensuring the absence of stagnation points  is recast in terms of the streamfunction $\psi$ as
\begin{equation}\label{umax_psi}
	\psi_y<0.
\end{equation}
The governing equations \eqref{euler_n} can then be reformulated (cf. \cite{Cons_book,CS2,CS11,Hen13,DH1} for details) in terms of the stream function as follows:
\begin{subequations}\label{Stream}
	\begin{align}\label{eqnsfirst}
		\Delta\psi &= \gamma(\psi)  & \mbox{ in } & -d<y<\eta(x),\\
 		\lvert \nabla \psi \rvert ^2+2g(y+d)&=Q  & \mbox{ on } & y=\eta(x), 
		\label{eulerorig}  
		\\ 
		\label{psi3}
		\psi&=0 & \mbox{ on } & y=\eta(x),	
		\\ 
		\label{stokeslast}
		\psi&=-p_0 & \mbox{ on } & y=-d,    
	\end{align}  
\end{subequations}
where $\gamma(\psi)=\omega$ is the vorticity function, $Q$ is a physical constant known as the total hydraulic head, and  $p_0$ is a constant of motion known as the relative mass-flux, defined by
\[
	p_0=\int_{-d}^{\eta(x)}(u(x,y)-c)dy<0.
\] 
The appropriate weak  formulation of the stream function governing equations \eqref{Stream} is then: 
\begin{subequations}\label{stream_weak}
\begin{align}
	\label{wstr1}
	(\psi_x\,\psi_y)_x-{1\over 2}(\psi_x^2-\psi_y^2)_y-\big(\tilde\Gamma
	(\psi/p_0)\big)_y&=0  & \mbox{ in } & D_{\eta},\\ 
	\label{wstr2}
	\psi &=-p_0  & \mbox{ on } & y=-d,  \\
	\label{wstr3} 
	\psi &=0 & \mbox{ on } & y=\eta(x),  \\
	\label{wstr4}
	\lvert \nabla \psi \rvert ^2+2g(y+d)&=Q & \mbox{ on } & y=\eta(x),
\end{align}
\end{subequations}  
where 
	$\tilde\Gamma(p)=\int_{0}^{p}p_0\gamma(s)ds$.
In this paper we are interested in solutions of \eqref{stream_weak} with
$\psi\in W^{2,r}_{per}(D_{\eta})\subset C^{1,\alpha}_{per}(\overline{D_{\eta}})$, $\tilde\Gamma\in C^{0,\alpha}
([-1,0])$, and $\eta\in C^{1,\alpha}_{per}(\R)$, for some $\alpha\in 
(0,1)$ and $r=2/(1-\alpha)$, where the boundary conditions \eqref
{wstr2}--\eqref{wstr4}  are satisfied in the classical sense, and  the 
equation \eqref{wstr1} is satisfied in the sense of distributions, and  $\psi_x,
\psi_y$ are also understood in the classical sense.

The third governing equation formulation we introduce is the form most suitable for local bifurcation theory to be implemented in order to prove the existence of water waves on a fluid domain of fixed mean-depth with a discontinuous vorticity distribution. It was first introduced in \cite{Hen13,DH1} to prove the existence of classical solutions representing small and large amplitude water waves propagating on a fluid of fixed mean-depth, where it was shown (cf.  \cite{Hen13,DH1}) that the governing equations \eqref{Stream} transform to
\begin{subequations}\label{H_0}
\begin{eqnarray}\label{hfirst}
\left(\frac{1}{d^2}+h^2_q\right) h_{pp}-2h_q(h_p+1)h_{pq}+(h_p+1)^2 h_{qq} 
+ \frac{\gamma(p)}{p_0}(h_p+1)^3 =0\nonumber \\ \mbox{ in } -1<p<0, \\ 
\label{hsecond}
\frac 1{d^2}+h_q^2+\frac{(h_p+1)^2}{p_0^2}[2gd(h+1)-Q]=0 \quad  \mbox{ on } p=0,
\\ \label{hbnd}
h=0 \quad  \mbox{ on }  p=-1.
\end{eqnarray} 
\end{subequations}
Here the variables $(q,p)$ are defined by the semi-hodograph transformation of variables
\begin{equation}\label{hod}
(x,y) \mapsto(q,p):= (x,\psi (x,y)/p_0),
\end{equation}
which has the effect of transforming the fluid domain with unknown free boundary $\eta$ into the fixed semi-infinite rectangular strip $\mathbb R\times 
[-1,0]$. It is now apparent that the non-stagnation condition  \eqref{umax} is vital in order to ensure that the change of variables \eqref{hod} represents an isomorphism.
 The function $h$, defined  by
\begin{equation}\label{hdef}
h(q,p)=\frac{y}d-p,
\end{equation}
is a modified-height function  where it is understood that $y=y(q,p)$ is a function of the new $(q,p)-$variables. In this paper we seek solutions $h$ of \eqref{H_0} which are  even and $2\pi-$periodic in $q$, and the analogous  condition to \eqref{eta} which the modified-height function should satisfy is given by 
\begin{equation}\label{hint}
\int_{-\pi}^{\pi} h(q,0)dq=0.
\end{equation} The condition \eqref{umax} excluding stagnation points   becomes 
\begin{equation}\label{hmax}
h_p+1>0,
\end{equation} and consequently the system \eqref{H_0} is a uniformly elliptic quasilinear partial differential equation with oblique nonlinear boundary conditions. The appropriate weak form of the modified-height function formulation \eqref{H_0} is given in divergence form by
\begin{subequations}\label{H}
\begin{align}\label{hfirst_0}
\left\{-\frac{1+d^2h^2_q}{2d^2(1+h_p)^2}+\frac{\Gamma(p)}{2d^2}\right\}_p+
\left\{\frac{h_q}{1+h_p} \right\}_q&=0  &\mbox{ in }  -1&<p<0, \\ \label{hsecond_0}
-\frac{1+d^2h^2_q}{2d^2(1+h_p)^2} -\frac{gd(h+1)}{p_0^2}+\frac{Q}
{2p_0^2}&=0 & \mbox{ on }  p&=0,
\\ \label{hbnd_0}
h&=0 &  \mbox{ on }  p&=-1,
\end{align} 
\end{subequations}
where 
\begin{equation*}\label{gamma}
\Gamma(p)=2\int_0^p\frac{d^2\gamma(s)}{p_0}ds, \quad -1\leq p \leq 0,
\end{equation*}
and in the following $\Gamma_{min}:= \min\limits_{p\in [-1,0]}\Gamma(p)$. We understand by a solution of \eqref{H} a function $h\in W^{2,r}_{per}(R)\subset C^{1,\alpha}_{per}(\overline R)$, where $r>2/(1-\alpha)$ for $\alpha\in(0,1)$ and $R=[-\pi,\pi]\times [-1,0]$. Here the \emph{per} subscript indicates that our solutions are even and $2\pi-$periodic in the $q-$variable, and the inclusion relation (along with the choice of $r$) are derived from Morrey's inequality.
\subsection{Equivalence of formulations}
In \cite{DH1} it was proven that the Euler equation \eqref{euler_n}, the stream function \eqref{Stream} and the modified-height function \eqref{H_0} formulations of the governing equations are equivalent in the sense of  classical  solutions. This approach can be modified (in the same vein as  in \cite{CS11} for the standard height function formulation) to prove the equivalence of the weak versions \eqref{euler_weak}, \eqref{stream_weak} and \eqref{H} of these formulations that we are interested in. Furthermore, we note that results concerning the equivalence of the standard height formulation of the  governing equations were recently strengthened  in \cite{VZ12}, and this approach was suitably adapted in  \cite{SSeq16} for the modified height function form of our equations giving the following result:
\begin{theorem}\cite{SSeq16} \label{equiv:form}
	Let $1/3<\alpha< 1$. Then the following 	formulations are equivalent:
	\begin{enumerate}
	\item[(i)] the weak velocity formulation \eqref{euler_weak} with  
	\eqref{umax}, for  $\eta\in C^{1,\alpha}_{per}(\R)$, and 
	$u,v, P\in C^	{0,\alpha}_{per}(\overline{D_{\eta}})$;
	\item[(ii)]  the weak stream function formulation \eqref{stream_weak} with 
	\eqref{umax_psi}, 
	for $\tilde{\Gamma}\in C^{0,\alpha}([-1,0])$,  $\eta
	\in C^{1,\alpha}_{per}(\R)$ and  $\psi\in C^{1,\alpha}_{per}(
	\overline{D_{\eta}})$;  
	\item[(iii)]  the weak height formulation \eqref{H} 
	together with \eqref{hmax}, for $\Gamma\in C^{0,\alpha}([-1,0])$ and  
	$h\in C^{1,\alpha}_{per}\big(\overline{R}\big)$.
	\end{enumerate}
\end{theorem}
\subsubsection{Regularity properties of solutions}
Recently, instigated by the paper \cite{ConEschAn}, a wide body of research  has flourished in establishing strong {\em a priori} regularity properties for periodic  water waves with vorticity in a variety of physical settings, such as flows with stratification \cite{HMreg}, free-boundary problems incorporating the capillarity effects of surface tension \cite{DHreg}, and flows with very general vorticity distributions \cite{EschMatAn,MMexist}. We note in particular that the approach and results of \cite{EschMatAn}, which apply to vorticity functions which are merely integrable, may be adapted for our modified height function formulation \eqref{H}, giving the following result: 
\begin{theorem}[Regularity of streamlines]
Given a solution $h \in C^{1,\alpha}_{per}\big(\overline{R}\big)$ of \eqref{H} such that \eqref{hmax} holds, then: 
\begin{itemize}
\item $ \partial_q^m h \in C^{1,\alpha}_{per}\big(\overline{R}\big)$ for all $m\in \mathbb N$.
\item Furthermore, there exists a constant $L>1$ such that, for all integers $M\geq 3$, we have
\[
\|\partial_q^m h\|_{1,\alpha}\leq L^{m-2}(m-3)!.
\]
\item Accordingly, all the streamlines for solutions of the Euler equation formulation \eqref{euler_weak}, including the free-surface wave profile $\eta(x)$, are {\em a priori} real analytic curves.
\end{itemize}
\end{theorem}
\subsection{Main Result}
The principle aim of this paper is to use local bifurcation theory to prove the following:
\begin{theorem}\label{MR}
	Let $0<\alpha<1$ and $\Gamma\in C^{0,\alpha}\left([-1,0]\right)$, and suppose we are given the wave speed $c>0$ and a fixed mean-water depth $d>0$. Suppose for $\lambda^*>-\Gamma_{min}$ there exists a weak solution $0\not\equiv M(p)\in C^{1,\alpha}([-1,0])$ to the Sturm-Liouville problem
	\begin{equation}\label{MR-SL}
\left\{
\begin{array}{ll}
	\displaystyle
	(a^3 M_{p})_p=d^2aM					& \mbox{ in } 	-1<p<0
	\\
	\displaystyle a^3M_p={gd^3\over p_0^2}M 		&  	\mbox{ on } p=0
	\\
	M=0 										&	\mbox{ on } p=-1,
\end{array}
\right.
\end{equation} where $a(\lambda,p)=\sqrt{\Gamma(p)+\lambda}$ for $p\in[-1,0]$. Then for  sufficiently small $\eps>0$ there exists a $C^1-$curve 
	$\mathcal{C}_{loc}=\{(\lambda^s,h^s)\in\R\times C^{1,\alpha}_{per}(\overline{R}): |s|<\eps\}$ of solutions to the water wave problem \eqref{H}, subject to $h_p+1>0$. The solution curve $\mathcal{C}_{loc}$ contains precisely one function independent of $q$, namely $\left(\lambda^s(0),h^s(0)\right)=\left(\lambda^*,H(p,\lambda^*)\right)$ where the laminar flow solution $H(p,\lambda^*)$ is defined by \eqref{laminar:sol} below. For $|s|<\eps$ we have 
\begin{equation}\label{Sol} h^s(q,p)=H(p,\lambda^s)+sM(p)\cos(q)+o(s)\quad \mbox{ in } C^{1,\alpha}(\overline{R}),
\end{equation} where $M$ is the solution to \eqref{MR-SL} with $M(0)=1$, and $Q$ is specified in terms of $\lambda$ by \eqref{Q:of:lambda}.
\end{theorem}
\noindent In Section \ref{SecLB} below we demonstrate how the existence of non-trivial water wave solutions of the form \eqref{Sol} is determined by system \eqref{MR-SL}. In the process, it will be clear that the existence of solutions to \eqref{MR-SL} depends on the nature of the vorticity distribution, something we expand further upon in Section \ref{SecEx}. In Theorem \ref{Gexist} it is shown that solutions of \eqref{MR-SL} exist for a wide range of general vorticity distributions, yet at the same time there also exist vorticity distributions for which solutions of \eqref{MR-SL} do not exist. In particular, we present necessary and sufficient conditions  for solutions of \eqref{MR-SL} to exist for constant, and some piecewise constant, vorticity distributions. 

We further note that the ``physical'' solutions $u,v,\eta$ of the Euler equations formulation \eqref{euler_weak} which correspond to \eqref{Sol} are periodic in the $x-$variable, with a single crest and trough per period, with $u,\eta$ being symmetric (and $v$ anti-symmetric) functions about the crestline. The inherent symmetry of these waves is not particularly prohibitive since its has been rigorously proven under quite general conditions that  water waves with vorticity are intrinsically symmetric, cf. \cite{Cons_book,ConEhWah,ConEschSym,MMsym}.

\section{Local bifurcation}\label{SecLB}
\subsection{The Crandall-Rabinowitz local bifurcation theorem}
The primary tool we  employ to prove existence of small amplitude water waves for the system \eqref{H}  is the Crandall-Rabinowitz  \cite{CranRab} local bifurcation 
theorem, which may be stated as follows (cf. \cite{Buff} for a detailed proof and discussion regarding the theorem and its hypotheses): 
\begin{theorem}[Crandall-Rabinowitz]\label{CR}
Let $X,Y$ be Banach spaces and let $\mathcal F\in C^k( \mathbb R \times X,Y)$ 
with $k\geq2$ satisfy:
\begin{enumerate}
\item \label{cr1} $\mathcal F(\lambda,0)=0$ for all $\lambda\in\mathbb R$;
\item \label{cr2} The Fr\'echet derivative $\mathcal F_x(\lambda^*,0)$ is a 
Fredholm operator of index zero with a one-dimensional kernel: \[
\ker (\mathcal F_x(\lambda^*,0))=\{sx_0: s\in\mathbb R, 0\neq x_0\in X\};\]
\item \label{cr3} The tranversality condition holds: 
\[
\mathcal F_{\lambda x}(\lambda^*,0)(1,x_0)\not \in \mathcal{R}(\mathcal 
F_x(\lambda^*,0)).
\]
\end{enumerate}
Then $\lambda^*$ is a bifurcation point in the sense that there exists $
\epsilon_0>0$ and a branch of solutions 
\[
\{
(\lambda,x)=\{(\Lambda(s),s\chi(s)): s\in \mathbb R, |s|<\epsilon_0\} \subset X 
\times \mathbb R,
\},
\]
with $\mathcal F(\lambda,x)=0$, $\Lambda(0)=\lambda^*,\chi(0)=x_0$, and the 
maps 
\[
s\mapsto \Lambda(s)\in \mathbb R,\qquad  s\mapsto s\chi(s)\in X,
\] are of class $C^{k-1}$ on $(-\epsilon_0,\epsilon_0)$.
Furthermore there exists an open set $U_0\in  \mathbb R \times X$ with 
$(\lambda^*,0)\in U_0$ and 
\[
\{
(\lambda,x)\in U_0:\mathcal F(\lambda,x)=0, x\neq 0\}=\{(\Lambda(s),s\chi(s)): 0<|
s|<\epsilon_0
\}.
\]
\end{theorem}

\subsection{Operator reformulation of system \eqref{H}}\label{Lam}
\noindent
The first step in our implementation of Crandall-Rabinowitz' theory is to determine solutions $H(p)$ of system \eqref{H} which have no $q-$dependence, and which therefore describe parallel shear flows with a flat horizontal surface and streamlines. In the process  we will elicit a suitable bifurcation parameter, we will establish a reformulation of system \eqref{H} as an operator between appropriate Banach spaces, and we will achieve a characterisation of the ``trivial'' solutions which fulfil the first condition of Theorem \ref{CR}. We solve \eqref{H} along similar lines as in \cite{DH1} to get, for $-\Gamma_{\min}<\lambda<Q$,  a $q-$independent solution of the form
\begin{align}
\nonumber H(p;\lambda)=\int_0^p\frac {ds}{\sqrt{\lambda+\Gamma(s)}}+\frac{1}{2gd}\left[Q-\frac{p_0^2}{d^2}\lambda\right]-(p+1), \quad -1<p\leq 0,
\\
\label{laminar:sol}
	=\int_{-1}^p{1\over \sqrt{\lambda+\Gamma(s)}}ds-(p+1)\in 
	C^{1, \alpha}_{per}\left([-1,0]\right).
\end{align}
The parameter $\lambda$ is related to the fluid velocity at the flat surface by the relation
\Beq \label{disp1}
\sqrt \lambda=\left.\frac{1}{H_p+1}\right|_{p=0}=\left.\frac{d(u-c)}{p_0}\right|_{on\ the\ flat\ surface},
\Eeq 
and it is implicitly related to $Q$ by the following formula
\begin{eqnarray}\label{Q:of:lambda}
Q=2gd\int_{-1}^0\frac {ds}{\sqrt{\lambda+\Gamma(s)}}+\frac{p_0^2}{d^2}\lambda>0.
\end{eqnarray}
We see that $Q(\lambda)$ is a positive convex function of $\lambda$ 
with a minimum occurring at the unique value $\lambda_0>0$ satisfying
\begin{equation}\label{Qmin}
{p_0\over gd^3}=\int_{-1}^0 {1\over \big(\lambda_0+\Gamma(s)\big)^{3\over 2}}ds,
\end{equation}
and $Q(\lambda)$ is monotonically decreasing for $-\Gamma_{min}<
\lambda<\lambda_0$ and monotonically increasing for $\lambda>\lambda_0$. We infer that the appropriate bifurcation parameter for system \eqref{H} is $\lambda>-\Gamma_{min}$. To re-express  system \eqref{H} as an  operator $\mathcal F(\lambda,x):\mathbb R\times X \rightarrow Y$, where  
$X,Y$ are Banach spaces, we work as follows. Let us consider first the operators 
\[
\begin{array}{ll}
	\mathcal{G}_1 h		&	\displaystyle=\left\{-\frac{1+d^2h^2_q}
	{2d^2(1+h_p)^2}+\frac{\Gamma(p)}	{2d^2}\right\}_p+\left\{\frac{h_q}
	{1+h_p} \right\}_q,  
	\\ 
	\mathcal{G}_2 h	&	\displaystyle=\left.\left(-\frac{1+d^2h^2_q}
	{2d^2(1+h_p)^2} -\frac{gd(h+1)}{p_0^2}+\frac{Q}{2p_0^2}\right)\right|_T
\end{array}
\]
where $T=\{(q,0): q\in [-\pi,\pi]\}$. With $\alpha\in(0,1)$ fixed we define the Banach spaces 
\[
	X=\left\{h\in C^{1,\alpha}_{per}(\overline{R}): \int_T h(q,0)dq=0, \, h=0 \mbox{ on }-1,\, h \mbox
	{ is even in }q\right\},
\]
and $Y=Y_1\times Y_2$, with 
\[
\begin{array}{l}
	Y_1=\left\{f\in \mathcal{D}'_{per}(R): \, f=\partial_q\varphi_1+\partial_p\varphi_2, 
	\mbox{ for some } \varphi_1,\varphi_2\in C^{0,\alpha}_{per}(\overline{R}), f 
	\mbox{ even in } q\right\}
	\smallskip\\
	Y_2=\left\{f\in C^{0,\alpha}_{per}(T): \, f \mbox{ is even in } q\right\}
\end{array}	
\]
where $\mathcal{D}'_{per}$ are the distributions which 
are $2\pi$-periodic in the $q$-variable and the norm in $Y_1$ is given by 
\[
	\|f\|_{Y_1}=\inf\limits_{\varphi_1,\varphi_2}\big\{\|\varphi_1\|_{C^{0,\alpha}_{per}(\overline{R})}
	+\|\varphi_2\|_{C^{0,\alpha}_{per}(\overline{R})}\big\}.
\]
Now, we define the operators $\mathcal{F}_i:\R\times X\to Y_i$, for 
$i=1,2$, by
\begin{equation}\label{Op}
	\mathcal{F}_1(\lambda,w)=\mathcal{G}_1(H+w)
	\qquad \mathcal{F}_2(\lambda,w)=\mathcal{G}_2(H+w),
\end{equation}
and the operator associated to \eqref{H} is  $\mathcal{F}:=(\mathcal{F}
_1,\mathcal{F}_2):\R\times X\to Y$. Since $\mathcal F(\lambda,0)=0$ for all $\lambda\in\R$ by design, the first hypothesis in the Crandall-Rabinowitz Theorem \ref{CR} concerning the existence of trivial solutions for the operator \eqref{Op} is fulfilled. 


\subsection{The linearised problem}\label{The linearised problem}

The linearisation of problem \eqref{H} is achieved by considering solutions of the form 
$h(q,p)=H(p,\lambda)+\eps m(q,p)$, where (at the first order of $\epsilon$) $m\in C^{1,\alpha}_{per}(\overline{R})$ must satisfy the problem 
\begin{equation}\label{lin:prob:1}
\left\{
\begin{array}{ll}
	\displaystyle
	(a^3 m_{p})_p+d^2a\,m_{qq}=0 				& 	\mbox{ in } R,
	\\
	\displaystyle a^3m_p={gd^3\over p_0^2}m 		&  	\mbox{ on } p=0,
	\\
	m=0 										& 	\mbox{ on } p=-1,
\end{array}
\right.
\end{equation}
with $m$ even and $2\pi$-periodic in the $q$-variable, and $a(\lambda,p)={1\over H_p+1}=\sqrt{\Gamma(p)+\lambda}$~ for $\lambda>-\Gamma_{min}$, cf. \cite{DH1} for details. We are interested in solutions of the type
$$
    m(q,p)=M(p)\cos(kq),
$$ 
 with $k\ge 1$ an integer and $M(p)\in C^{1,\alpha}\left([-1,0]\right)$, and we infer that $M$ should solve the Sturm-Liouville problem 
\begin{equation}\label{Sturm:Liouv}
\left\{
\begin{array}{ll}
	\displaystyle
	(a^3 M_{p})_p=k^2d^2aM 					& \mbox{ in } 	-1<p<0,
	\\
	\displaystyle a^3M_p={gd^3\over p_0^2}M 		&  	\mbox{ on } p=0,
	\\
	M=0 										&	\mbox{ on } p=-1.
\end{array}
\right.
\end{equation}
Problem \eqref{Sturm:Liouv} has an associated variational formulation which can be expressed in terms of 
\[
	\mathbb{F}(\lambda,\varphi)={\displaystyle -gd^3\varphi^2(0)+p_0^2\int_{-1}^0 
	a^3\varphi_p^2 dp\over \displaystyle p_0^2d^2\int_{-1}^0a\varphi^2 dp},
\]
whereby the minimization problem associated to \eqref{Sturm:Liouv} is given by 
\begin{equation}\label{mu:lambda} 
	\mu(\lambda)=
	\inf\limits_{{\varphi\in H^1\left(-1,0\right),\atop \varphi(-1)=0,\,\varphi\not
	\equiv 0}}\mathbb{F}(\lambda,\varphi).
\end{equation}
The function $M$ attaining the minimum \eqref{mu:lambda} is smooth and 
satisfies the Sturm-Liouville problem
\[
	(a^3M_p)_p=-\mu d^2aM\quad \mbox{in } (-1,0),
\]
with boundary conditions given by \eqref{Sturm:Liouv}, cf. \cite{DH1} for details. Since  
$a$ and $\mu$ are $C^1$-functions of $\lambda$, the following results transfer over from the setting of classical solutions considered in \cite{DH1} unimpeded:
\begin{lemma}\cite{DH1}\label{monotonicity:mu} {\bf(Monotone)}
	 $\mu(\lambda)$ is a strictly increasing function of $\lambda$ when
	 $\mu(\lambda)<0$.
\end{lemma}
\noindent Furthermore, given any vorticity distribution $\gamma$  it can be shown using a reasoning similar to Section 5.23 of \cite{DH1} that  $\mu(\lambda)>-1$ for sufficiently large $\lambda$. Piecing this information together with Lemma \ref{monotonicity:mu} gives us the following result:
\begin{corollary}\label{lambda_unique_sol}
	A neccessary and sufficient condition for the existence of a solution of \eqref{lin:prob:1} is that $\mu(\lambda)\leq -1$ for some $\lambda>-\Gamma_{min}$. Furthermore, a solution of $\mu(\lambda^*)=-1$, for $\lambda^*>-\Gamma_{min}$, must also be unique. 
\end{corollary}
\noindent The following result (which follows along the same lines as Section 6.1 of \cite{DH1}) will be useful in later considerations:
\begin{lemma}\cite{DH1}\label{L0}
Let $\lambda=\lambda_0$ be the unique value where $Q(\lambda)$ (defined by \eqref{Q:of:lambda}) attains its minimum. Then $\mu(\lambda_0)=0$ and furthermore relation \eqref{Qmin}, which may be expressed as
\[
\int_{-1}^0a^{-3}(\lambda,s)ds=\frac{p_0^2}{gd^3},
\] 
holds uniquely for $\lambda=\lambda_0$. 
\end{lemma}
\noindent In order to ascertain when the final two hypotheses in the Crandall-Rabinowitz Theorem \ref{CR} hold, we first calculate $\mathcal{F}_{w}$ by linearising the operator $\mathcal{F}$ around $w=0$: 
\begin{subequations}\label{LinF}
\begin{equation}\label{linearised:op:1}
\begin{array}{ll}
	\displaystyle \left.\mathcal{F}_{1,w}\right|_{w=0}
	&\displaystyle=\left\{ -{d^22w_q\partial_q 2d^2
	(1+H_p+w_p)^2-4d^2(1+H_p+w_p)\partial_p
	(1+d^2w_q^2) \over (2d^2)^2(1+H_p+w_p)^4}\right\}_p
	\\
	&\displaystyle
	\left.+\left\{{\partial_q(1+H_p+w_p)-w_q\partial_p\over 
	(1+H_p+w_p)^2}\right\}_q\right|_{w=0}
	\smallskip\\
	&\displaystyle= \left\{{\partial_p\over d^2(1+H_p)^3}\right\}_p+\left\{
	{\partial_q\over 1+H_p}\right\}_q
\end{array}
\end{equation} 
and 
\begin{equation}\label{linearised:op:2}
\begin{array}{ll}
	\displaystyle \left.\mathcal{F}_{2,w}\right|_{w=0}
	&\!\!\!\!\displaystyle=\left.\left.\left( -{ 2d^2w_q\partial_q 2d^2
	(1+H_p+w_p)^2-4d^2(1+H_p+w_p)(1+d^2w_q^2) 
	\partial_p \over (2d^2)^2(1+H_p+w_p)^4}-{gd\over p_0^2}Id\right)
	\right|_{w=0}\right|_T
	\smallskip\\
	&\!\!\!\!\displaystyle= \left.\left({\partial_p\over d^2(1+H_p)^3}-{gd\over 
	p_0^2}I\right)\right|_T.
\end{array}
\end{equation} 
\end{subequations}
\subsection{Null space of $\mathcal{F}_{w}(\lambda^*,0)$}\label{SubsecNull}
The aim of  this section is to establish that, when a solution to \eqref{MR-SL} exists, the null space of $\mathcal F_{w}(\lambda^*,0)$, denoted by $\ker (\mathcal F_{w}(\lambda^*,0))$, is one-dimensional. This is the first step in establishing that the second hypothesis in Crandall-Rabinowitz' Theorem \ref{CR} is satisfied.   
\begin{lemma}\label{Nullspace}{\bf (Nullspace)} 
	If $-\Gamma_{min}<\lambda^*<\lambda_0$ is the unique solution of $\mu(\lambda)
	=-1$, then the null space of $\mathcal{F}_{w}(\lambda^*,0)$ is 
	one-dimensional.
\end{lemma}
\begin{proof}
	From Section \ref{The linearised problem}  we know that the solution 
	 $m(q,p)=M(p)\cos(q)$ of \eqref{lin:prob:1}, where $M$ is 
	the eigenfunction corresponding to the eigenvalue $k=1$ of  
	\eqref{Sturm:Liouv} with $M(0)=1$, is contained in the null space of $\mathcal{F}_{w}(\lambda^*,0)$. 
	To see that $\ker\left(\mathcal{F}_{w}(\lambda^*,0)\right)$ is generated solely by $M(p)\cos(q)$ we let $m\in \ker\left(\mathcal{F}_{w}(\lambda^*,0)\right)$, then for every fixed $p\in [-1,0]$ we write 	\[
		m(q,p)=\sum\limits_{k=0}^{\infty} m_{k}(p)\cos(kq) \mbox{ in } L^2\left([-
		\pi,\pi]\right),
	\]  
	where the coefficients  $m_k\in C^{1,\alpha}\left([-1,0]\right)$ are given by 
	\[
		m_0(p)={1\over 2\pi}\int_{-\pi}^{\pi}m(q,p)dq,\qquad m_{k}(p)=
		{1\over \pi}\int_{-\pi}^{\pi}m(q,p)\cos(kq)dq,\; k\ge 1.
	\]
	Since $m\in \ker\left(\mathcal{F}_{w}
	(\lambda^*,0)\right)$ it follows that $m$ satisfies \eqref{lin:prob:1}, and multiplying \eqref{lin:prob:1}
	by $\cos(kq)$ and integrating over $[-\pi,\pi]$ we find that each coefficient $m_k$ 
	must satisfy \eqref{Sturm:Liouv}. 	
For $k\ge 2$ we must have $m_k\equiv 0$,  since otherwise if $m_k\not\equiv 0$ for some $k\ge 2$  the minimising property of $\mu(\lambda^*)=-1$ is contradicted by the relation
	\[
		\mathbb{F}(\lambda^*,m_k)={-gd^3m_k(0)+p_0^2\int_{-1}^0a^3(\partial_p 
		m_k)^2dp\over p_0^2d^2\int_{-1}^0 am_k^2 dp}=-k^2<-1.
	\]
        Regarding $m_0$, since it solves \eqref{Sturm:Liouv} 
	with $k=0$ it must satisfy the differential equation 
	$$
			[a^3 (m_0)_p]_p=0,
	$$
	and integrating over $(-1,p)$ we infer that 
	there exists some constant $C_0\in \R$ such that 
	\[
		m_0(p)=C_0\int_{-1}^p a^{-3}ds.
	\]
	The boundary condition in \eqref{Sturm:Liouv} at $p=0$ gives us
	\[
		C_0=C_0{gd^3\over p_0^2}\int_{-1}^0 a^{-3}ds,
	\]
which, considering Lemma \ref{monotonicity:mu}, Lemma \ref{L0} and relation \eqref{Qmin}, can only hold if $C_0=0$ and accordingly $m_0=0$.
  Hence we are left with the function $m_1(p)$, which must be  a constant multiple of $M(p)$, and so
	$$m=m_1(p)\cos(q)=C M(p)\cos(q), \mbox{ where } C\in \R,$$ 
	and the null space of 	$\mathcal{F}_{w}(\lambda^*,0)$ 
	is indeed one-dimensional, generated by the eigenfunction corresponding 
	to the eigenvalue $k=1$ of \eqref{Sturm:Liouv}. 
	\end{proof}

\subsection{The range of $\mathcal{F}_{w}(\lambda^*,0)$}\label{SubsecRange}
In this section we give a characterization of the range of $\mathcal{F}
_{w}(\lambda^*,0)$, where $-\Gamma_{min}<\lambda^*<\lambda_0$ satisfies $\mu(\lambda^*)=-1$. The proof of Proposition \ref{range} exhibits slight differences from similar scenarios presented in \cite{CS11,DH1}, yet these subtle differences are vitally important, indeed crucial, in enabling the finely balanced apparatus of local bifurcation theory to be applied to the operator \eqref{Op}. Accordingly,  for completeness we present the proof in its entirety.
\begin{proposition}\label{range}{\bf (Range)}
	The pair $(\mathcal{A},\mathcal{B})\in Y_1\times Y_2$, with $\mathcal{A}=\partial_q \mathcal{A}_1 +\partial_p 
	\mathcal{A}_2$ belongs to the range of $\mathcal{F}_{w}(\lambda^*,0)$ if and only if 
	\begin{equation}\label{rel:A:B}
		-\sum\limits_{j=1}^2 \int\int_R \mathcal{A}_j\partial_j\varphi^* dq dp=\int_{T}
		\varphi^* \mathcal{B} dq
	\end{equation}
	where $\varphi^*\in X$ 
	generates the null space of $\mathcal{F}_{w}(\lambda^*,0)$.
\end{proposition}
\begin{proof} {\em Necessity:}
	let $(\mathcal{A},\mathcal{B})$ belong to the range of $\mathcal{F}_{w}(\lambda^*,0)$, 
	then there exists $\Phi\in X$ such that 
	\begin{equation}\label{range:A:B}
	\left\{
		\begin{array}{ll}
			\displaystyle \mathcal{A}=\left\{{\Phi_p\over d^2(1+H_p)^3}\right\}_p+
			\left\{{\Phi_q\over (1+H_p)}\right\}_q & \mbox{in }R,
			\\
			\displaystyle \mathcal{B}={\Phi_p\over d^2(1+H_p)^3}-{gd\over p_0^2}\Phi  
			& \mbox{on }	T.
		\end{array}
	\right.
	\end{equation}
	Multiplying $\mathcal{A}$ by $\varphi^*\in  \ker (\mathcal{F}_
	{w}(\lambda^*,0))$, integrating over $R$, then it follows immediately (since $\Phi$ and $\varphi^*$ are even in $q$) that \eqref{rel:A:B} holds.

\noindent {\em Sufficiency:} Given $\left(\mathcal{A,B}\right)\in Y_1 \times Y_2$ such that \eqref{rel:A:B} holds, we must prove that there exists a $\Phi\in X$ such that \eqref{range:A:B} is satisfied. A vital component of the proof of existence will hinge on an application of the Lax-Milgram Theorem \cite{Ev}, and due to subtle coercivity considerations (similar to those outlined in \cite[p.61]{Cons_book}) we decompose our analysis in terms of the zero, and non-zero, Fourier modes separately, working as follows. 
    
\noindent {\em Zero Fourier mode}: Given $(\mathcal{A},\mathcal{B})$, since $\int_{-\pi}^{\pi}\varphi^*dq=0$ it is trivial to see that the zero Fourier modes
	\[
		\mathcal{B}_0={1\over 2\pi}\int_{-\pi}^{\pi} \mathcal{B}(q)dq,\quad 
		\mathcal{A}_0(p)={1\over 2\pi}
		\partial_p\int_{-\pi}^{\pi} \mathcal{A}_2(p,q)dq,	\] 
	satisfy relation \eqref{rel:A:B} in their own right. We observe that if relation \eqref{range:A:B} were to hold for $\mathcal{A,B}$, then integrating \eqref{range:A:B} with respect to $q$ from $-\pi$ to $\pi$ we find that $(\mathcal{A}_0,\mathcal{B}_0)$ and  
	$\Phi_0={1\over 2\pi}\int_{-\pi}^\pi \Phi(q,p)dq$ should satisfy  
	\begin{equation}\label{sol:A:0:B:0}
	\left\{
	\begin{array}{ll}
		\displaystyle \mathcal{A}_0=\left\{{\partial_p\Phi_{0}\over d^2(1+H_p)^3}\right\}_p
		&	\mbox{in }(-1,0),
		\smallskip\\
		\displaystyle
		\mathcal{B}_0={\partial_p\Phi_{0}\over d^2(1+H_p)^3}-{gd\over p_0^2}\Phi_0
		&	\mbox{on }p=0.
	\end{array}
	\right.
	\end{equation}
	With this in mind, we prove that equation  \eqref{sol:A:0:B:0} has a unique solution $\Phi_0$, thereby settling the matter of sufficiency for the zero Fourier modes of the range. Integrating 
	\eqref{sol:A:0:B:0} over $(-1,p)$ we have
	\begin{equation}\label{sol:A:0:B:0:a:0}
	\begin{array}{rcl}
		\displaystyle \int_{-1}^p \mathcal{A}_0(\tau)d\tau	&	=	&	\displaystyle
		{\Phi'_{0}(p)\over d^2(1+H_p(p))^3}-	{\Phi'_{0}(-1)\over d^2(1+H_p
		(-1))^3},
			\end{array}	
	\end{equation}
		which upon further integration over $(-1,p)$ gives us 		
	\begin{equation}\label{sol:A:0:B:0:a}
	\Phi_0(p)={C}\int_{-1}^p(1+H_p^3(s))ds+\int_{-1}^p d^2
		(1+H_p(s))^3\int_{-1}^s \mathcal{A}_0(\tau)d\tau ds,			\end{equation}
	for a constant $C={\Phi'_{0}(-1)/(1+H_p(-1))^3}$. Considering  \eqref{sol:A:0:B:0:a:0} and \eqref{sol:A:0:B:0:a} 
	at 	$p=0$ gives  
	\begin{equation*}
	\begin{array}{rcl}
		\displaystyle \int_{-1}^0 \mathcal{A}_0(\tau)d\tau	&	=	&	\displaystyle
		{\Phi'_{0}(0)\over d^2(1+H_p(0))^3}-{C\over d^2},

\smallskip \\		\displaystyle \int_{-1}^0 d^2(1+H_p(s))^3\int_{-1}^s \mathcal{A}_0(\tau)d\tau ds
		&	=	&	
		\displaystyle \Phi_0(0)-{C}\int_{-1}^0(1+H_p(s))^3ds,
		\end{array}	
	\end{equation*}
respectively, and plugging the above relations into the boundary condition at $p=0$ in \eqref{sol:A:0:B:0}  gives 
	\begin{equation*}
	\begin{array}{l}
		\displaystyle C\left({1\over d^2}- {gd\over p_0^2}\int_{-1}^0(1+H_p)^3
		ds\right)
		=	
		\displaystyle	\mathcal{B}_0-\int_{-1}^0 \!\!\!\mathcal{A}(\tau)d\tau
		+{gd^3\over p_0^2}\int_{-1}^0 d^2(1+H_p)^3\int_{-1}^s 
		\mathcal{A}_0(\tau)d\tau ds.
	\end{array}
	\end{equation*}
	Lemma \ref{L0}, together with relation \eqref{Qmin},  ensures that the terms in parentheses on the left-hand side above are non-zero, and therefore we can determine the constant $C$ uniquely. Accordingly, problem \eqref{sol:A:0:B:0} has a unique solution determined by \eqref{sol:A:0:B:0:a}. 	

\noindent {\em Non-zero Fourier mode}: For the non-zero Fourier modes, we note that given a solution $\Phi$ of \eqref{range:A:B}, and a corresponding function $\Phi_0$ which solves \eqref{sol:A:0:B:0}, the function $v=\Phi-\Phi_0\in X_0$ solves  
	\begin{equation}\label{range:A:B_X_0}
	\left\{
		\begin{array}{rl}
			\displaystyle \left\{{v_p\over d^2(1+H_p)^3}\right\}_p+
			\left\{{v_q\over (1+H_p)}\right\}_q=\mathcal{A}-\mathcal{A}_0 & \mbox{in }R,
			\\
\displaystyle
			{v_p\over d^2(1+H_p)^3}-{gd\over p_0^2}v=\mathcal{B}-\mathcal{B}_0  
			& \mbox{on }	T, 
		\end{array}
	\right.
	\end{equation}
where the space $X_0$ is defined by
\[ X_0=\left\{\varphi\in X:\, \int_{-\pi}^\pi \varphi (q,p)dq=0 \,
	\mbox{ for all }\; p\in [-1,0]\right\},\]
and  $(\mathcal{A}-\mathcal{A}_0, \mathcal{B}-\mathcal{B}_0)\in Y_0$ where 
\[ Y_0	=\left\{(\mathcal{A},\mathcal{B})\in Y: \, \partial_p
	\int_{-\pi}^\pi \mathcal{A}_2(q,p)
	dq=0\;\mbox{ for all }\; p\in
	[-1,0],\;\int_{-\pi}^\pi \mathcal{B}(q)dq=0\right\}.
\]
Therefore, bearing in mind the solution \eqref{sol:A:0:B:0:a}, the question of proving the existence of a solution to the full problem \eqref{range:A:B} can be reduced instead to that of solving the system \eqref{range:A:B_X_0}. We proceed by first seeking solutions $v^{(\eps)}\in X_0$ of the auxiliary equations  
	\begin{equation}\label{approx:eq}
	\left\{
		\begin{array}{ll}
			\displaystyle \mathcal{A}=-\eps v^{(\eps)}+(1+\eps)\left\{{v^{(\eps)}_p
			\over 
			d^2
			(1+H_p)^3}\right\}_p+(1+\eps)
			\left\{{v^{(\eps)}_q\over (1+H_p)}\right\}_q & \mbox{in }R,
			\\
			\displaystyle \mathcal{B}=(1+\eps){v^{(\eps)}_p\over d^2(1+H_p)^3}-{gd
			\over p_0^2}v^{(\eps)}  
			& \mbox{on }	p=0,
			\\
			v^{(\eps)}=0 & \mbox{on }	p=-1,
			\displaystyle
		\end{array}
	\right.
	\end{equation}
	and we observe that if  $v^{(\eps)}\in W_{per}^{1,2}(\overline{R})$ is a weak solution
	of \eqref{approx:eq} then the relation
	\begin{equation}\label{WEAK}
		\begin{array}{l}
		\displaystyle\int\!\!\!\int_R \left[(1+\eps){v^{(\eps)}_q\over 1+H_p}
		\varphi_q+(1+\eps){v^	{(\eps)}_p\over d^2(1+H_p)^3}\varphi_p+\eps 
		v^{(\eps)}\varphi\right]dq\,dp-\int_{T}{gd\over p_0^2}v^{(\eps)}\varphi dq
		\smallskip\\
		\displaystyle = \int\!\!\!\int_R (\mathcal{A}_1\varphi_q+\mathcal{A}_2\varphi_p)dq\,dp+\int_T (\mathcal{B}-
		\mathcal{A}_2)\varphi dq
		\end{array}
	\end{equation}
	must hold for all test functions $\varphi\in C^{1}_{per}(\overline{R})$ which are zero 
	on $p=-1$, where $\mathcal{A}=\partial_q \mathcal{A}_1+\partial_p \mathcal{A}_2$ for $\mathcal{A}_1,\mathcal{A}_2\in 
	C^{0,\alpha}_{per}(\overline{R})$. 
 Defining the Hilbert space
	\[
		\mathbb{H}=\left\{\varphi\in W^{1,2}_{per}(\overline{R}): \, \int_{-\pi}^\pi \varphi(q,p)
		dq=0 \mbox{ a.e. on } [-1,0], \varphi=0 \mbox{ a.e. on }p=-1\right\},
	\]
	then similar to the reasoning employed in \cite{DH1} it can be demonstrated that the left-hand side of \eqref{WEAK} defines a bounded, coercive bilinear form on $\mathbb{H}$--- indeed, coercivity considerations are the primary motivation for transposing our analysis of system \eqref{range:A:B} on the space $X$, to focussing instead on system \eqref{range:A:B_X_0} applied to the subspace $X_0$. Furthermore, as in \cite{DH1} it follows that the right-hand side of \eqref{WEAK} defines a bounded linear functional on 
	$\mathbb{H}$. 
Hence an implementation of the Lax-Milgram Theorem ensures the existence of a unique weak solution 
	$v^{(\eps)}\in \mathbb{H}$ to the approximate equations \eqref{approx:eq}, and following elliptic regularity considerations (as in, for example, Section 4.2 in \cite{CS11}) we may infer additional regularity in the solution, namely $v^{(\eps)}\in C^{1,\alpha}_{per}(\overline{R})$, and additionally an associated Schauder estimate holds,
	\begin{equation}\label{estimate:v:eps}
	\|v^{(\eps)}\|_{ C^{1,\alpha}_{per}(\overline{R})}\le C\{\|\mathcal{A}\|_
		{Y_1}+\|\mathcal{B}\|_{Y_2}+\|v^{(\eps)}\|_{L^{\infty}(R)}\},
	\end{equation} where the constant $C$ depends  only on $\|H_p+1\|_{C^{0,\alpha}_{per}(\overline R)}$.
	We will conclude our proof of the existence of a solution to the system \eqref{range:A:B_X_0} by proving that there exists a subsequence $\{\epsilon_n\}\rightarrow 0$ whose associated solutions $\{v^{(\eps_n)}\}_{n\ge 1}$ of \eqref{approx:eq} converge to some function $v$.
First we show that for any 	$\eps_n\to 0$, the sequence $\big\{\|v^{(\eps_n)}\|_{C^{1,\alpha}_{per}	(\overline{R})}\big\}_{n\ge 1}$ is bounded. 
	Arguing  by contradiction, if 
	$\|v^{(\eps_n)}\|_{C^{1,\alpha}_{per}(\overline{R})}\to \infty$ 
	as $\eps_n\to 0$ then \eqref{estimate:v:eps} yields $\|v^{(\eps_n)}\|_
	{L^{\infty}(R)}\to\infty$, while the functions $u_n={v^
	{(\eps_n)}\over \|v^{(\eps_n)}\|_{L^{\infty}(R)}}$ are bounded in $C^{1,
	\alpha}_{per}(\overline{R})$. Accordingly, due to the compact embedding $ C^{1,\alpha}_
	{per}(\overline{R})\subset\subset  C^{1}_{per}(\overline{R})$ we can  
	extract a subsequence $\{u_{n_k}\}$ that converges  in $C^1_{per}
	(\overline{R})$ to some $u$ with $\|u\|_{L^{\infty}(R)}=1$  as $n_k\to \infty$. Upon dividing \eqref{approx:eq} by $\|v^{(\eps_{n_k})}\|_
	{L^{\infty}(R)}$ and taking limits we see that $u$ must solve (weakly) 
	\begin{equation}\label{prob:limit:u}
	\left\{
		\begin{array}{ll}
			\displaystyle 0=\left\{{u_p
			\over d^2(1+H_p)^3}\right\}_p+
			\left\{{u_q\over (1+H_p)}\right\}_q & \mbox{in }R,
			\\
			\displaystyle 0={u_p\over d^2(1+H_p)
			^3}-	{gd\over p_0^2}u  
			& \mbox{on }	p=0,
			\\
			u=0 & \mbox{on }	p=-1,
			\displaystyle
		\end{array}
	\right.
	\end{equation}
	and due to the Schauder estimate \eqref{estimate:v:eps}  we have an immediate gain of regularity, namely $u\in C^{1,\alpha}(\overline{R})$. Comparing \eqref{prob:limit:u} with \eqref{LinF} we infer that $u\in \ker \mathcal{F}_{w}(\lambda^*,0)$, and Lemma 
	\ref{Nullspace} implies that $u=\delta \varphi^*$ for some $\delta\in\R$. Choosing $\varphi=\varphi^*$ in \eqref{WEAK} and bearing in mind \eqref
	{rel:A:B} we obtain 
	\begin{equation}\label{weak:form:v:eps:n}
		(1+\eps_n)\int\!\!\!\int_{R}\left[{\varphi^*_q\over 1+H_p}v_q^
		{\eps_n}+
		{\varphi^*_p\over d^2(1+H_p)^3}v_p^{(\eps_n)}\right]dq\,dp+\eps_n
		\int\!\!\!\int_R v^{(\eps_n)}\varphi^* dq dp={gd\over p_0^2}\int_Tv^
		{(\eps_n)}\varphi^*  dq,
	\end{equation}
	which following integration by parts and using the fact that $\varphi^*$ 
	solves \eqref{prob:limit:u} gives 
	\begin{equation}\label{weak:form:v:eps:n:2}
		\eps_n\int\!\!\!\int_R v^{(\eps_n)}\varphi^* dq dp+\eps_n{gd\over 
		p_0^2}\int_T v^{(\eps_n)}\varphi^* dq=0.
	\end{equation}
	Dividing by $\eps_n\|v^{(\eps_n)}\|_{L^{\infty}(R)}$, and taking limits  gives 
	\[
		\delta\int\!\!\!\int_R \big(\varphi^*\big)^2 dq dp+\delta{gd\over 
		p_0^2}\int_T \big(\varphi^*\big)^2 dq=0,
	\] since $u_n={v^
	{(\eps_n)}/ \|v^{(\eps_n)}\|_{L^{\infty}(R)}}\to\delta\varphi^*$, which implies that $\delta=0$ and hence  $u=\delta\varphi^*=0$. This contradiction of the identity $\|u\|_{L^{\infty}(R)}=1$ leads us to conclude that the sequence $\big\{\|v^
	{(\eps_n)}\|_{C^{1,\alpha}_{per}(\overline{R})}\big\}_{n\ge 1}$ is indeed bounded as $\epsilon_n \rightarrow 0$, which together with the compact embedding 
	$C^{1,\alpha}_{per}(\overline{R})\subset
	\subset C^{1}_{per}(\overline{R})$ enables us to infer the existence of a subsequence 
	of $\{v^{(\eps_n)}\}_{n\ge 1}$ converging in $C^{1}_{per}(\overline
	{R})$ to some $v\in C^{1}_{per}(\overline{R})$ which is a solution to 
	\eqref{range:A:B_X_0}. By standard elliptic regularity \cite
	{Gilbarg} we have that $v\in C^{1,\alpha}_{per}(\overline{R})
	$. 

\noindent We have therefore proven that, when \eqref{rel:A:B} holds, there exists unique solutions $\Phi_0$ and $v$ of \eqref{sol:A:0:B:0} and \eqref{range:A:B_X_0} respectively, and $\Phi=\Phi_0+v$ then is a unique solution of \eqref{range:A:B}. This proves the sufficiency of relation \eqref{rel:A:B}. 
\end{proof}
The considerations of this section, coupled with those of subsection \ref{SubsecNull}, prove that when $\mu(\lambda^*)=-1$ the second hypothesis of Theorem \ref{CR} is satisfied.
%

\subsection{Transversality condition}
 To complete the proof that local bifurcation occurs for the simple eigenvalue $\lambda^*$, thereby proving our main result Theorem \ref{MR}, we must prove that  the transversality condition, hypothesis three of Theorem \ref{CR}, holds whenever $\mu(\lambda^*)=-1$ .
We recall that the transversality condition states that 
\[
	\mathcal{F}_{w,\lambda}(\lambda^*,0)\varphi^*\not\in \mathcal{R}
	\left(\mathcal{F}_{w}(\lambda^*,0)\right),
\]
where $\mathcal{F}_{w,\lambda}(\lambda^*,0)$ is given by
\[
	\mathcal{F}_{w,\lambda}(\lambda^*,0)w=\left({3\over 2d^2}(a
	w_p)_p+{1\over 2}(a^{-1}w_q)_q,\left.{3\over 2d^2}a
	w_p\right|_T\right),
\]
for $a=\sqrt{\Gamma(p)+\lambda}$. As a result of Proposition \ref{range} we only need to check that 
\begin{equation}\label{transv:cond}
	{1\over 2}\int\!\!\!\int_R \left[\big(a^{-1}\varphi^*_q\big)_q+
	{3\over d^2}\big(a\varphi^*_p\big)_p\right]\varphi^*dq\,dp-{3\over 2 
	d^2}\int_T a\varphi^*_p\varphi^*dq\ne 0.
\end{equation}
Integrating \eqref{transv:cond} by parts, and using the fact that $\varphi^*=0$ on $p=-1$, 
we have
\[
	-{1\over 2}\int\!\!\!\int_R a^{-1}\big(\varphi^*_q\big)^2 dq\, dp-
	{3\over d^2}\int\!\!\!\int_Ra\big(\varphi^*_p\big)^2dq\,dp<0.
\] 
Thus the  transversality condition is satisfied, and this completes the proof of our main result.

\section{Applications of Theorem \ref{MR-SL}}\label{SecEx}
In this section we examine the applications of our local existence result Theorem \ref{MR} with regard to  vorticity distributions, both general and specific. In particular, for arbitrary  vorticity distributions we present conditions which are sufficient for Theorem \ref{MR} to apply, whereas for some examples of constant, and piecewise-constant, vorticity distributions we obtain physically explicit conditions which possess the additional value of  being necessary for local bifurcation to occur. 
  \subsection{General vorticity distributions}
We recap that it was shown in \cite{DH1} for the setting of classical solutions of the water wave problem (where the vorticity distributions are at least continuous) that a sufficient condition for bifurcation to occur  is given by the following: 
\begin{proposition}\cite{DH1}
Suppose that $\gamma$ is an (otherwise arbitrary) vorticity distribution which satisfies
\Beq
\label{condition} \frac{\sqrt 2}{3}\gamma_{\infty}^{\frac 32}|p_0|^{\frac{1}{2}}|p_1|^{\frac{1}{2}}+\frac{2\sqrt 2}{5}\gamma_{\infty}^{\frac 12}|p_0|^{\frac{3}{2}}|p_1|^{\frac{3}{2}}<g,
\Eeq
where $\gamma_{\infty}=\| \gamma \|_{C[-1,0]}$ and $p_1=\max\{p\in[-1,0]:\Gamma(p)=\Gamma_{min}\}$. Then local bifurcation occurs and we obtain non-trivial solutions to the classical water wave problem. 
\end{proposition}
\noindent The proof of the above theorem can be modified to give the following result for discontinuous vorticity distributions, which  we present a brief sketch of here for the sake of completeness.
\begin{theorem}\label{Gexist}
	For $0<\alpha<1$ and $\Gamma\in C^{0,\alpha}\left([-1,0]\right)$, let $\theta=
	\sup\limits_{p\ne p_1}{|\Gamma(p)-\Gamma(p_1)|\over |p-p_1|^{\alpha}}$. If
	\begin{equation}\label{cond:g}
		g>{\theta^{3/2}\over 6\alpha d^3}|p_0|^2|p_1|^{3\alpha/2-1}+
		{\theta^{1/2}\over (2+\alpha/2)d}|p_0|^2|p_1|^{\alpha/2+1},	
	\end{equation}
	 then Theorem \ref{MR} holds and hence local bifurcation occurs.
\end{theorem}
\begin{remark}\label{REM}
We note that $\gamma \in L^{\infty}[-1,0]$ implies $\Gamma \in C^{0,1}[-1,0]$ is Lipschitz continuous, and if $\gamma(p)\geq0$ then $p_1=0$ by definition. Hence, a consequence of Theorem \ref{Gexist} is that for any  positive, bounded, but otherwise arbitrary vorticity distribution, local bifurcation  always occurs. We note that this remark applies trivially to irrotational flow, for which $\gamma \equiv 0$.
\end{remark}
\begin{proof}
We first note that  since
	\[
		0\le \Gamma(p)-\Gamma_{min}=\Gamma(p)-\Gamma(p_1)\le |p-p_1|^{\alpha}\cdot \sup\limits_{p\ne p_1}{\Gamma(p)-\Gamma(p_1)\over |p-p_1|^
		{\alpha}}= \theta 
		|p-p_1|^{\alpha},
	\] 
	we have that 
	\[
		a(-\Gamma_{min},p)=\sqrt{\Gamma(p)-\Gamma_{min}}\le \sqrt
		{\theta|p-p_1|^{\alpha}},\quad p\in [-1,0].
	\]
	From \eqref{cond:g}, we can take $k>1/2$ so that 
	\begin{equation}\label{cond:g_k}
		{gd^3 \over p_0^2}>\theta^{3/2}{1\over 2k-1+3\alpha/2}|p_1|^{3\alpha/
		2-1}+d^2\theta^{1/2}{1\over 2k+\alpha/2+1}|p_1|^{\alpha/2+1}.	
	\end{equation} 
	We define 
	\begin{equation}\label{def_p_n}
		p_n=\left({1-{1\over n}}\right)p_1-{1\over n}\in [-1,p_1],\; n\ge 2
	\end{equation}
	and  
	\begin{equation}\label{}
		\varphi_n(p)=\left\{
		\begin{array}{ll}
			0,			&	-1\le p\le p_n\\
			(p-p_n)^k,	&	p_n\le p\le 0.
		\end{array}
		\right.
	\end{equation}
	for all $n\ge 2$. It follows that $-1\le p\le p_n\le p_1\le 0$, and since 
	$k>1/2$ we have that $\varphi_n\in H^1\left([-1,0]\right)$. From the 
	definition of $p_n$ we have that $\varphi_n(-1)=0$ and  
	$\varphi_n\not\equiv 0$. 
	By way of direct computations we obtain 
	\begin{align*}
		\displaystyle\int_{-1}^0 a^3(-\Gamma_{min},p)\big
		(\partial_p 
		\varphi_n(p)\big)^2dp+&\int_{-1}	^0 d^2a(-\Gamma_{min},p)\varphi_n^2
		(p)dp
		\smallskip\\
		\displaystyle &\le \theta^{3/2}k^2\left[{|p_n|^{2k-1+3\alpha/
		2}\over 2k-1+3\alpha/2}+{3\alpha/2 |p_1-p_n|^{2k-1+3\alpha/
		2}\over (2k-1)(2k-1+3\alpha/2)}\right]
		\smallskip\\
		\displaystyle &
		+d^2\theta^{1/2}\left[{|p_n|^{2k+1+\alpha/
		2}\over 2k+1+\alpha/2}+{\alpha/2|p_1-p_n|^{2k+1+\alpha/
		2}\over (2k+1)(2k+1+\alpha/2)}\right]
		\smallskip\\
		\displaystyle & = |p_n|^{2k}\big(A_n+B_n\big)=\varphi_n^2(0)\big(A_n
		+B_n\big),
	\end{align*}
	where 
	\[
	\begin{array}{l}
		\displaystyle A_n=\theta^{3/2}k^2{|p_n|^{3\alpha/
		2-1}\over 2k-1+3\alpha/2}+d^2\theta^{1/2}{|p_n|^{\alpha/
		2+1}\over 2k+1+\alpha/2},
		\\
		\displaystyle B_n=\theta^{3/2}k^2{3\alpha/2 |p_1-p_n|^
		{3\alpha/
		2-1}\over (2k-1)(2k-1+3\alpha/2)}+d^2\theta^{1/2}{\alpha/2|p_1-p_n|^
		{\alpha/2+1}\over (2k+1)(2k+1+\alpha/2)}.
	\end{array}
	\]
	From the definition  \eqref{def_p_n} of $p_n$ we have that $|p_1-p_n|\le |p_n|$, and $\lim
	\limits_{n\to\infty}|p_1-p_n|=0$. From \eqref{cond:g_k} 
	there exists $N\ge 2$
	and $\eps>0$ such that $A_n<gd^3/p_0^2-\eps$ for all $n\ge N$. With this specified 
	value of $\eps$, we choose $n\ge N$ large enough such that $B_n<
	\eps$. With the choice of this $n$,  the function $\varphi_n$ satisfies
	\[
		\int_{-1}^0a^3(-\Gamma_{min},p)\big
		(\partial_p 
		\varphi_n(p)\big)^2dp+\int_{-1}	^0 d^2a(-\Gamma_{min},p)\varphi_n^2
		(p)dp<\frac{gd^3}{p_0^2}\varphi_n^2(0).
	\]
	Therefore, there exists $\lambda>-\Gamma_{min}$ satisfying $\mu(\lambda)
	<-1$, hence Corollary \ref{lambda_unique_sol} ensures the existence of a solution of \eqref{MR-SL}, which in turn through Theorem \ref{MR} implies that local bifurcation occurs.
\end{proof}

\subsection{Constant vorticity}
We observe that a very convenient implication of Remark \ref{REM} is that local bifurcation occurs for any positive constant vorticity $\gamma\geq0$. In \cite {DH1}  it was demonstrated that, for  constant, continuous negative vorticity $\gamma<0$, bifurcation occurs if and only if 
\begin{equation}\label{XX1}
\gamma^2d^2<(g+\gamma^2d)\tanh(d).
\end{equation}
\subsection{Piecewise-constant vorticity}
We present two examples here of different physical scenarios which are modelled by piecewise continuous vorticity. The first case, comprises a uniform layer of constant non-zero vorticity located at the surface of the flow, and beneath this layer the flow is purely irrotational. This model has an interesting physical relevance, since is well known that wind-generated water waves exhibit a localised layer of vorticity near the surface. In the initial stages of wind generated-waves on the ocean, the wind ``grips'' the surface, producing small-amplitude water waves which subsequently develop in size.   The second case is, in a sense, the opposite set-up, consisting of a layer of constant non-zero vorticity at the bottom of the fluid, and above this layer there is no vorticity. As a physical motivation for this type of flow,  it is well-known that wave-current interaction (and hence vorticity) is most dramatically observed where this is a rapid change in the current strength, and the near-bed region is a location which commonly experiences currents, accounting for sediment transport, for instance.
\noindent We note that, as a direct consequence of Remark \ref{REM}, local bifurcation occurs in both examples if $\gamma \geq 0$.  In \cite{HDisp13,HDisp} the Sturm-Liouville problem \eqref{MR-SL} was analysed for both models with the principle purpose being to obtain dispersion relations for the flow by way of a detailed analysis of relation \eqref{disp1}. The dispersion relations for small-amplitude waves detail how the relative speed of the wave at the free-surface varies with respect to certain physical flow parameters, such as the fixed mean-depth of the fluid, the wavelength, the vorticity distribution, and--- for the discontinuous vorticity distributions --- the location $\md$ of the isolated layer of vorticity.  Analogous results concerning  dispersion relations, obtained instead through analysing the standard height function formulation, were derived in \cite{ConDisp,CS11,dispM2,dispM1}. Aside from obtaining dispersion relations, a serendipitous boon of the approach employed in \cite{HDisp13,HDisp} is the derivation, when the region of vorticity is negative, of necessary and sufficient conditions for the existence of solutions to \eqref{MR-SL}. As a result of Theorem \ref{MR}, these double--up as necessary and sufficient conditions for local bifurcation to occur, and we outline them now below. 
\subsubsection{Case (a): Localised vorticity layer at the surface}
The vorticity distribution  takes the form
\Beqn 
\omega=\left\{   
\begin{array}{ll}
\gamma, & -\md <y<0, \\
0, & -d<y< -\md,
\end{array}
\right.
\Eeqn
where $\md \in [0,d]$ is some given fixed depth. It is understood by $\md=0$ in this context that the flow is entirely irrotational, and $\md=d$ signifies that the flow has a  constant, continuous vorticity distribution $\gamma$ throughout the entire domain. It was shown in \cite{HDisp13} that, for $\gamma<0$, a necessary and sufficient condition for both solutions of \eqref{MR-SL} exist and, consequent to Theorem \ref{MR}, that local bifurcation occurs, is given by
\Beq\label{k1}
\gamma^2{\md}\left(\md -\tanh(\md)\right)<g\tanh(\md).
\Eeq
We note that, upon taking the limit $\md \rightarrow d$, this condition for bifurcation particularises to the equivalent condition \eqref{XX1} for flows of constant, continuous negative vorticity.

\subsubsection{Case (b): Localised vorticity layer at the bed}
 For this case, the vorticity distribution 
\Beqn
\omega=\left\{   
\begin{array}{ll}
0, & -\md <y<0, \\
\gamma, & -d<y< -\md.
\end{array}
\right.
\Eeqn
In this setting, $\md=0$ represents purely rotational flow while for $\md=d$ the flow is irrotational. For $\gamma<0$, it was shown in \cite{HDisp} that a necessary and sufficient condition that there be a solution of \eqref{MR-SL}, and so local bifurcation occurs, takes the form  
\[
|\gamma|< \frac{\sqrt{g(\sinh\left(d\right)(d-\md) -\sinh\left(\md\right)\sinh\left(d-\md\right))}}{(d-\md)\sqrt{(d-\md)\cosh\left(d\right) -\gamma^2\cosh\left(\md\right)\sinh\left(d-\md\right)}}.
\]
While being altogether more intricate than \eqref{k1}, we note that upon taking the limit $\md \rightarrow 0$ this condition also reduces to the equivalent condition \eqref{XX1} for flows of constant and continuous negative vorticity.



\end{document}